\documentclass[11pt]{amsart}
\usepackage[margin=1in]{geometry}

\usepackage{amssymb}
\usepackage{amsthm}
\usepackage{amsmath}
\usepackage{mathrsfs}
\usepackage{amsbsy}
\usepackage[all]{xy}
\usepackage{bm}
\usepackage{hyperref}
\usepackage{tikz}
\usepackage{array}
\usepackage{float}
\usepackage{enumerate}
\usepackage{xcolor}
\usepackage{hhline}
\allowdisplaybreaks
\usepackage[noadjust]{cite}

\usepackage{caption}
\usepackage{subcaption}
\usepackage{diagbox}
\usepackage{booktabs}

\DeclareFontFamily{U}{rcjhbltx}{}
\DeclareFontShape{U}{rcjhbltx}{m}{n}{<->rcjhbltx}{}
\DeclareSymbolFont{hebrewletters}{U}{rcjhbltx}{m}{n}

\let\aleph\relax\let\beth\relax

\DeclareMathSymbol{\aleph}{\mathord}{hebrewletters}{39}
\DeclareMathSymbol{\beth}{\mathord}{hebrewletters}{98}

\usepackage[noabbrev,capitalise]{cleveref}

\newenvironment{enumerate*}%
  {\begin{enumerate}[(I)]%
    \setlength{\itemsep}{10pt}%
    \setlength{\parskip}{0pt}}%
  {\end{enumerate}}

\newtheorem{theorem}{Theorem}[section]
\newtheorem{proposition}[theorem]{Proposition}

\newtheorem{lemma}[theorem]{Lemma}

\theoremstyle{definition}

\title{The Ulam sequence grows subexponentially}

\author{Noah Kravitz}
\address{St John's College, Oxford and Mathematical Institute, University of Oxford; St Giles', Oxford OX1 3JP, UK}
\email{noah.kravitz@maths.ox.ac.uk}

\begin{document}

\maketitle

\begin{abstract}
We prove that the $n$-th element of the Ulam sequence, a somewhat mysterious recursively-defined sequence of positive integers, has size $e^{O(\sqrt{n})}$; previous upper bounds were exponential in $n$.  Our argument uses an additive-combinatorial lemma, possibly of independent interest, about covering finite difference sets with multiplicatively short intervals.  The proofs are due to ChatGPT.

\end{abstract}

\section{Introduction}
In 1964, Stanislaw Ulam~\cite{ulam} introduced the following curious integer sequence that has since come to be known as the \emph{Ulam sequence}.  Start with $u_1:=1$ and $u_2:=2$; then, for each $n \geq 3$, let $u_n$ be the smallest integer greater than $u_{n-1}$ that has a unique representation as a sum of two \emph{distinct} previous elements of the sequence.  The first few terms of the Ulam sequence are $$1,2,3,4,6,8,11,13,16,18,26.$$  Write $\mathcal{U}:=\{u_1, u_2, u_3, \ldots \}$ for the set of all of the elements of the Ulam sequence.

Despite the simplicity of Ulam's definition, the behavior of the sequence $u_n$ is somewhat erratic and mysterious.  Even the asymptotic growth rate remains unknown.  Computations in the 1960's initially suggested that $u_n$ might grow super-linearly in $n$, but more extensive recent computations indicate that a linear growth rate is more likely.  For more on this problem's surprisingly extensive historical background, see the discussions in~\cite{stefan,CS} (see also~\cite{me,hin,hs,andrei} for some generalizations to settings other than the integers).

The trivial upper bound, from the observation $u_n \leq u_{n-1}+u_{n-2}$, is that $u_n$ is at most the $(n+1)$-th Fibonacci number.  This is an exponential upper bound, where the base of the exponential is the golden ratio (around $1.618$).  Successive improvements reduced the base of the exponential, with the record of $1.454$ due to Cl\'ement and Steinerberger in 2025~\cite{CS}.  ChatGPT (see AI statement below) has provided a short and clever argument for a subexponential upper bound.

\begin{theorem}\label{thm:main}
There is a constant $C>0$ such that $u_n \leq e^{C \sqrt{n}}$ for all $n$.
\end{theorem}

The proof uses only the following two properties of the Ulam sequence:
\begin{itemize}
    \item consecutive elements of $\mathcal{U}$ have quotient at most $2$;
    \item every element of $\mathbb{N} \setminus \mathcal{U}$ has either zero or at least two different representations as a sum of two distinct elements of $\mathcal{U}$.
\end{itemize}
In particular, the argument uses only the fact that the Ulam sequence includes every positive integer with a unique representation as a sum of two distinct smaller elements of the sequence, and it
does not make use of the fact that the Ulam sequence omits every positive integer with multiple representations as a sum of two distinct smaller elements of the sequence.  The set of positive integers with at most two $1$'s in their binary expansions also satisfies the two bullet-pointed conditions, and its $n$-th smallest element is $2^{\Theta(\sqrt{n})}$.  So to obtain a qualitatively stronger result, one would have to make essential use of the omission from $\mathcal{U}$ of integers with multiple representations in $\mathcal{U}+\mathcal{U}$.

The proof idea, in a nutshell, is to show that if $\mathcal{U}$ has many large gaps in the interval $[1,x]$, then it also contains many elements either moderately smaller than $x$ or slightly larger than $x$; in other words, if $\mathcal{U}$ has many sparse spots at a given scale, then it must have some dense spots nearby.  One then averages this statement over different choices of $x$.

\section{Additive-combinatorial input}
The proof of Theorem~\ref{thm:main} makes use of an independently interesting additive-combinatorial fact about covering finite difference sets with multiplicatively short intervals.

\begin{proposition}\label{prop:add}
There is a constant $C>0$ such that the following holds.  Let $0<\varepsilon\leq 0.99$, and let $A \subseteq \mathbb{R}$ be a finite set.  Then $(A-A) \cap \mathbb{R}_{>0}$ can be covered by at most $$C\varepsilon^{-1/2} (\log \varepsilon^{-1})^{1/2} |A|$$ intervals of the form $[z,(1+\varepsilon)z)$.
\end{proposition}

We remark that for our application to Theorem~\ref{thm:main}, it would suffice to have a bound of the form $\varepsilon^{-1+\gamma}|A|$ for any fixed $\gamma>0$.  The quantitative dependence in the conclusion of the proposition is tight up to a constant factor; see Appendix~\ref{apdx} for a matching construction.

One can interpret $(A-A) \cap \mathbb{R}_{>0}$ as the set of distances among the elements of $A$.  Proposition~\ref{prop:add} also applies, with essentially an identical proof, to the more general setting of metric spaces: The set of (nonzero) distances among the elements of a finite set $A$ in a metric space can be covered by $O(\varepsilon^{-1/2} (\log \varepsilon^{-1})^{1/2} |A|)$ intervals of the form $[z,(1+\varepsilon)z)$.

A precursor of Proposition~\ref{prop:add} appears in an influential 1995 paper of Callahan and Kosaraju~\cite{CK}, which implies the proposition with $\varepsilon^{-1/2} (\log \varepsilon^{-1})^{1/2}$ replaced by $\varepsilon^{-1}$.  Some thematically similar ideas appear in the work of Bourgain~\cite{bourgain} on the discretized sum-product problem (see also~\cite{gkz}).

\subsection{A simpler result}

Proposition~\ref{prop:add} is most interesting for $\varepsilon$ tending to $0$.  The range of constant $\varepsilon$ can be handled by a simpler argument.  Say that a finite set of real numbers is \emph{dissociated} if all of its subset sums are distinct.  The \emph{additive dimension} of a set is the largest size of a dissociated subset.  We observe that if $A,B \subseteq \mathbb{R}$ are finite sets, then $A-B$ has additive dimension at most $|A|+|B|-1$.  To see this, let $$\{a_1-b_1, \ldots, a_m-b_m\}$$ be a dissociated subset of $A-B$.  Consider the bipartite graph on the vertex set $A \sqcup B$, where $a_i,b_i$ are adjacent for each $1 \leq i \leq m$.  We claim that this graph is acyclic.  Indeed, any cycle would (after relabeling) correspond to a relation of the form
$$(a_1-b_1)-(a_2-b_2)+\cdots+(a_{2\ell-1}-b_{2\ell-1})-(a_{2\ell}-b_{2\ell})=0$$
for some $\ell \leq m/2$, but this is impossible by dissociativity.  Thus the graph is a forest, and its number $m$ of edges is at most $|A|+|B|-1$, as desired.  (A similar idea appears in~\cite{shk}.)

An immediate consequence of this observation is that $(A-B) \cap \mathbb{R}_{>0}$ can be covered by at most $|A|+|B|-1$ intervals of the form $[z,2z)$.  Indeed, in a left-to-right greedy disjoint covering by such intervals, taking the smallest element of $A-B$ in each interval gives a dissociated set, so the number of such intervals is at most $|A|+|B|-1$.  Moreover, each interval $[z,2z)$ can be covered by $O(\varepsilon^{-1})$ intervals of the form $[z,(1+\varepsilon)z)$, so more generally $(A-B) \cap \mathbb{R}_{>0}$ can be covered by $O(\varepsilon^{-1} (|A|+|B|))$ intervals of the form $[z,(1+\varepsilon)z)$.  The point of Proposition~\ref{prop:add} is improving the $\varepsilon$ dependence.  (To obtain an asymmetric statement about $A-B$, one can apply Proposition~\ref{prop:add} with $A \cup B$ playing the role of $A$.)

\subsection{Proof of the proposition}
We finally turn to the proof of Proposition~\ref{prop:add}.  The new character entering our story here is $K(t)$, the minimum number of intervals of the form $[z,z+t)$ needed to cover the set $A$.  We analyze $A-A$ at different scales one at a time.  Let $R:=100\varepsilon^{-1}$.  Suppose we want to understand the elements of $(A-A) \cap [R^i,R^{i+1})$ for some $i \in \mathbb{Z}$.  Such elements can be written as $a-b$, where $a,b \in A$ differ by around $R^i$.  Take a covering of $A$ by $K(R^{i-1})$ intervals of length $R^{i-1}$.  Then $a,b$ must lie in different intervals.  If instead we want to cover $A$ by intervals of length $R^{i+2}$, then we can use one long interval to cover both of the short intervals that originally covered $a,b$.  So elements of $A-A$ at scale $R^i$ force $K(R^{i+2})$ to be smaller than $K(R^{i-1})$.

We will see that even though a single ``merging'' of length-$R^{i-1}$ intervals into a length-$R^{i+2}$ interval can correspond to many elements of $(A-A) \cap [R^i,R^{i+1})$, it cannot correspond to too many multiplicatively $(1+\varepsilon)$-separated elements of $(A-A) \cap [R^i,R^{i+1})$.  Making this argument precise yields that every multiplicatively $(1+\varepsilon)$-separated subset of $(A-A) \cap [R^i,R^{i+1})$ has size $$\ll \varepsilon^{-1/2}(\log \varepsilon^{-1})^{1/2}(K(R^{i-1})-K(R^{i+2})).$$
Summing over $i \in \mathbb{Z}$ gives a telescoping sum that yields a further factor of $|A|$.  The details follow.

\begin{proof}[Proof of Proposition~\ref{prop:add}]
Define the scale $$R:=100\varepsilon^{-1}.$$
Let $i \in \mathbb{Z}$, and consider $(A-A) \cap [R^i, R^{i+1})$.  Extract a maximum-size multiplicatively $(1+\varepsilon)$-separated subset $X_i$ of $(A-A) \cap [R^i, R^{i+1})$ greedily from left to right.  Then $X_i$ is additively $\varepsilon R^i$-separated, and $(A-A) \cap [R^i, R^{i+1})$ is covered by the $|X_i|$ intervals of the form $[z,(1+\varepsilon)z)$ for $z \in X_i$.

Take a covering $$\{I_1, \ldots, I_{K(R^{i-1})}\}$$ of $A$ by intervals of length $R^{i-1}$, listed from left to right.  Each element of $X_i$ is contained in $I_j-I_k$ for some $j<k$; fix one such pair $(j,k)$ for each element of $X_i$.  Since $I_j-I_k$ is an interval of length $2R^{i-1}<\varepsilon R^{i}$, it contains at most one element of $X_i$.  So we have exactly $|X_i|$ distinct pairs $(j,k)$ associated to the elements of $X_i$.  Consider the graph $G$ on the vertex set $[K(R^{i-1})]$ with this edge set.  Within each connected component of $G$, consecutive intervals have their left endpoints at most $R^{i+1}+R^{i-1}<R^{i+2}-R^{i-1}$ apart.  So the intervals in a connected component of size $s \geq 2$ can be covered by at most $\lceil s/2 \rceil$ intervals of length $R^{i+2}$.  Note that $\lceil s/2 \rceil \leq s-s/3$ for $s \geq 2$.  Since $G$ has $|X_i|$ edges in total, its connected components of size at least $2$ must together span at least $|X_i|^{1/2}$ (say) vertices.  It follows that
$$K(R^{i+2}) \leq K(R^{i-1})-|X_i|^{1/2}/3.$$
At the same time, we trivially have $|X_i| \ll \varepsilon^{-1} \log \varepsilon^{-1}$ by multiplicative volume-packing.  Thus
$$|X_i| \ll \varepsilon^{-1/2}(\log \varepsilon^{-1})^{1/2}|X_i|^{1/2}\ll \varepsilon^{-1/2}(\log \varepsilon^{-1})^{1/2}(K(R^{i-1})-K(R^{i+2})).$$

Summing over $i \in \mathbb{Z}$, we find that $(A-A) \cap \mathbb{R}_{>0}$ can be covered with
$$\sum_{i \in \mathbb{Z}}|X_i| \ll \varepsilon^{-1/2}(\log \varepsilon^{-1})^{1/2} \sum_{i \in \mathbb{Z}}(K(R^{i-1})-K(R^{i+2}))$$
intervals of the form $[z,(1+\varepsilon)z)$.  Note that $K(R^i)=1$ for $i$ sufficiently large, and $K(R^i)=|A|$ for $i$ sufficiently small.  So the previous sum telescopes to
$$\sum_{i \in \mathbb{Z}}(K(R^{i-1})-K(R^{i+2}))=3(|A|-1) \ll |A|,$$
as desired.
\end{proof}

Note that our choice of $R$ was dictated by the constraint $2R^{i-1}<\varepsilon R^i$ in the second paragraph of the proof.  For a smaller choice of $R$, one would need to replace the ``short'' scale $R^{i-1}$ with $R^{i-h}$ for some $h$ large enough that $2R^{i-h}<\varepsilon R^i$.  Then our volume-packing upper bound would be $|X_i| \ll \varepsilon^{-1} \log R$, and the final telescoping step would pick up a factor $h+2$ instead of $3$.  Carrying through these changes would lead to a bound of the shape
$$\sum_{i \in \mathbb{Z}}|X_i| \ll \varepsilon^{-1/2} (\log R)^{1/2} \left(\frac{\log \varepsilon^{-1}}{\log R}+1\right)|A|,$$
where the two terms are balanced when $R=\varepsilon^{-\Theta(1)}$.

\section{Proof of Theorem~\ref{thm:main}}

Throughout this section, fix a large positive integer $n$.  We will upper-bound $u_n$.

\subsection{Large and small gaps}
We start from the telescoping identity
$$\log u_n=\sum_{i=2}^n \log(u_i/u_{i-1}).$$
It is natural to look at multiplicative telescoping (rather than additive telescoping) since we are working on a scale that is closer to exponential than to polynomial.  Note that every $u_i \leq 2u_{i-1}$.  Our main business will be showing that $u_i/u_{i-1}$ is rarely large.  To this end, define the sets
$$S(\delta):=\{u_i: (1-\delta)u_i \geq u_{i-1}, ~i \in [2,n]\}$$
for $0< \delta \leq 1/2$, which record which $u_i$'s are preceded by large gaps.  We now rewrite our telescoping identity as
\begin{align*}
\log u_n &=\sum_{i=2}^n \log(u_i/u_{i-1})=\sum_{i=2}^n \int_1^{u_i/u_{i-1}} \frac{1}{x} \,dx\\
 &=\int_1^2 \frac{|\{2 \leq i \leq n: u_i/u_{i-1} \geq x\}|}{x} \,dx\\
 &=\int_{0}^{1/2} \frac{|S(\delta)|}{1-\delta} \,d\delta \leq 2\int_{0}^{1/2} |S(\delta)| \,d\delta,
\end{align*}
where in passing to the third line we changed variables to $\delta:=1-1/x$.  We will show that
$$|S(\delta)|\ll \delta^{-3/4} (\log \delta^{-1})^{3/4} \cdot \sqrt{n}$$
(see Lemma~\ref{lem:s-delta-bound} below), whence we will conclude that
$$\log u_n \ll 2 \sqrt{n} \int_{0}^{1/2} \delta^{-3/4} (\log \delta^{-1})^{3/4} \,d\delta \ll \sqrt{n}$$
since the integral converges (for this convergence, we care only that the power of $\delta^{-1}$ is below $1$).

\subsection{The key combinatorial step}

Fix $0<\delta \leq 1/2$.  The crux is obtaining an upper bound on $|S(\delta) \cap (0,x)|$ in terms of the number of elements of $\mathcal{U}$ near $x$.  To do this, we show that many sums $x+y$ with $y \in S(\delta) \cap (0,\delta x)$ can \emph{fail} to be in $\mathcal{U}$ only if other elements of $\mathcal{U}$ cluster around $x$.
\begin{lemma}\label{lem:key}
Let $x \in S(\delta)$.  Then
$$|S(\delta) \cap (0,x)| \ll \delta^{-1/2} (\log \delta^{-1})^{1/2} \cdot |\mathcal{U} \cap (\delta x, x)|+\delta^{-1} \log \delta^{-1} \cdot |\mathcal{U} \cap [x,(1+\delta)x)|.$$
\end{lemma}

\begin{proof}
We can trivially bound
$$|S(\delta) \cap [\delta x,x)| \leq 1+\frac{\log \delta^{-1}}{\log ((1-\delta)^{-1})} \ll \delta^{-1} \log \delta^{-1},$$
since the elements of $S$ are separated by multiplicative factors of at least $(1-\delta)^{-1}$.  We now restrict our attention to the elements $y \in S(\delta) \cap (0,\delta x)$.

The number of such $y$'s with $x+y \in \mathcal{U}$ is at most $$|\mathcal{U} \cap (x,(1+\delta)x)|.$$  For the remaining $y$'s, there are two ways that $x+y$ could fail to be in $\mathcal{U}$.

First, consider the $y$'s where $x+y$ fails to be in $\mathcal{U}$ because $x+y$ has a second representation in $\mathcal{U}+\mathcal{U}$ with one summand greater than $x$. We obtain an expression
$$x+y=(x+c)+(y-c)$$
for some $0<c<y$ with $x+c, y-c \in \mathcal{U}$.  Note that in fact $$\delta y \leq c<y$$ since $y \in S(\delta)$ has a large gap before it.  The number of such $c$'s is trivially at most $|\mathcal{U} \cap (x,(1+\delta)x)|$ due to the constraints $x+c \in \mathcal{U}$ and $c<y<\delta x$.  Each $c$ corresponds to values of $y$ contained in the interval $(c,\delta^{-1}c]$, and the number of such $y$'s is $\ll \delta^{-1} \log \delta^{-1}$ by the multiplicative separation of $S(\delta)$.  Summing over the choices of $c$ gives a total contribution of $$\ll \delta^{-1} \log \delta^{-1} \cdot |\mathcal{U} \cap (x,(1+\delta)x)|.$$

Second, consider the $y$'s where $x+y$ fails to be in $\mathcal{U}$ because $x+y$ has a second representation in $\mathcal{U}+\mathcal{U}$ with both summands smaller than $x$.  We obtain an expression
$$x+y=a+b$$
for some distinct $a,b \in \mathcal{U}$ with $a,b<x$.  From $x \in S(\delta)$ we see that $a,b \leq (1-\delta)x$ and hence also $a,b> \delta x$ (since $a+b>x$).  So $y=a-(x-b)$ is contained in the difference set $A-A$ for
$$A:=(\mathcal{U} \cap (\delta x, x)) \cup (x-(\mathcal{U} \cap (\delta x, x))).$$
Proposition~\ref{prop:add} tells us that the set of such $y$'s can be covered by $$\ll \delta^{-1/2} (\log \delta^{-1})^{1/2} \cdot |\mathcal{U} \cap (\delta x, x)|$$ intervals of the form $[z,(1+\delta)z)$.  Each such interval contains at most a single $y$, again due to multiplicative separation.

We obtain the conclusion of the lemma by summing the final centered expressions from the previous four paragraphs and noting that $|\mathcal{U} \cap [x,(1+\delta)x)|=|\mathcal{U} \cap (x,(1+\delta)x)|+1$.
\end{proof}

One can understand this lemma conceptually as follows.  The element $x+y$ can fail to be in $\mathcal{U}$ only if it has another representation in $\mathcal{U}+\mathcal{U}$ as either $(x+d)+(y-d)$ or $(x-d)+(y+d)$ for some positive $d$.  The large gaps before $x,y$ (due to $x,y \in S$) force $d$ to be somewhat large.  This in turn restricts $y$ to a fairly small set, where ``fairly small'' depends on the elements of $\mathcal{U}$ near $x$.

\subsection{Clean-up}
Lemma~\ref{lem:key} gives information about $\mathcal{U}$ near a single element $x \in S(\delta)$.  Summing over $x \in S(\delta)$ and rearranging yields the desired upper bound on $|S(\delta)|$.

\begin{lemma}\label{lem:s-delta-bound}
We have $|S(\delta)|\ll \delta^{-3/4} (\log \delta^{-1})^{3/4} \cdot \sqrt{n}$.
\end{lemma}

\begin{proof}
Let $x_1<\cdots<x_R$ denote the elements of $S(\delta)$.  Apply Lemma~\ref{lem:key} for $x=x_j$ and sum over $1 \leq j \leq R-1$ to obtain
$$\binom{R-1}{2} \ll \delta^{-1/2} (\log \delta^{-1})^{1/2} \sum_{j=1}^{R-1} |\mathcal{U} \cap (\delta x_j, x_j)|+\delta^{-1} \log \delta^{-1} \sum_{j=1}^{R-1}|\mathcal{U} \cap [x_j,(1+\delta)x_j)|.$$
For the first sum, swapping the order of summation gives
$$\sum_{j=1}^{R-1} |\mathcal{U} \cap (\delta x_j, x_j)|=\sum_{k=1}^n |\{j \leq R-1: \,u_k<x_j<\delta^{-1}u_k\}|\ll n \delta^{-1} \log \delta^{-1},$$
as usual due to the multiplicative separation of $S(\delta)$.  For the second sum, notice that the intervals $[x_j,(1+\delta)x_j)$ are disjoint (since all elements of $S$ are preceded by long gaps) and contained in $[1,u_n]$.  Thus
$$\sum_{j=1}^{R-1} |\mathcal{U} \cap [x_j,(1+\delta)x_j)| \leq |\mathcal{U} \cap [1,u_n]|=n.$$
We conclude that
$$R^2 \ll 1+\binom{R-1}{2} \ll n \delta^{-3/2} (\log \delta^{-1})^{3/2},$$
and the result follows on taking square roots.
\end{proof}

We described at the beginning of this section how Theorem~\ref{thm:main} follows from Lemma~\ref{lem:s-delta-bound}.

\appendix

\section{Tightness of Proposition~\ref{prop:add}}\label{apdx}

In this appendix we demonstrate the tightness of Proposition~\ref{prop:add}.  For brevity, let $M_\varepsilon(A)$ denote the minimum number of intervals of the form $[z,(1+\varepsilon)z)$ needed to cover $(A-A) \cap \mathbb{R}_{>0}$.

We first record that if $A_1,A_2 \subseteq \mathbb{R}_{\geq 0}$ are finite sets satisfying $$\min\{(A_2-A_2) \cap \mathbb{R}_{>0}\} \geq 2 \max \{A_1\},$$ then $$M_{\varepsilon}(A_1 \cup A_2) \geq M_\varepsilon(A_1)+M_\varepsilon(A_2),$$ since no interval $[z,(1+\varepsilon)z)$ can simultaneously contain points of both $A_1-A_1, A_2-A_2$.  In particular, given $A \subseteq \mathbb{R}_{\geq 0}$, there is some $\lambda>1$ such that for all positive integers $n$, we have
$$\frac{M_\varepsilon(A \cup \lambda \cdot A \cup \cdots \cup \lambda^n \cdot A)}{|A \cup \lambda \cdot A \cup \cdots \cup \lambda^n \cdot A|} \geq \frac{M_\varepsilon(A)}{|A|}.$$
This shows that once we have a single example $A$ for which Proposition~\ref{prop:add} is tight up to constants, we can produce larger examples with the same property (for the same value of $\varepsilon$).
 
To see that $\varepsilon^{-1/2}$ is the necessary power of $\varepsilon$ in the proposition, one can take $N:=\lfloor \varepsilon^{-1/2}\rfloor$ and consider the $2N$-element set
$$A:=\{0,1, 2, \ldots, N-1\} \cup \{N, 2N, \ldots, N^2\},$$
which has $M_\varepsilon(A)=N^2\gg \varepsilon^{-1/2}|A|$.  Obtaining the further factor of $(\log \varepsilon^{-1})^{1/2}$ requires more care, as follows.

Let $N$ be a large positive integer.  Set $T:=N \cdot 8^N+1$, and define the set
$$A:=\{t/T: 1 \leq t \leq T\} \cup \{8^n(1+h/8^N): 1 \leq n \leq N, \, 1 \leq h \leq 8^N\}.$$
The first constituent part of $A$ is an arithmetic progression at a small scale, and the second constituent is essentially a perturbation of a geometric progression at a large scale.  Our set has size $$|A|=T+N \cdot 8^N \asymp N \cdot 8^{N}.$$

Consider the elements of $(A-A)\cap \mathbb{R}_{>0}$ of the form
$$x(n,h,t):=8^n(1+h/8^N)-t/T.$$
Note that $x(n,h,t) \in [(1/2) \cdot 8^n, 2 \cdot 8^n]$ subsists at scale $8^n$.  In particular, elements $x(n,h,t)$ with different values of $n$ are multiplicatively separated by factors of at least $2$, so they must be covered by disjoint sets of intervals $[z,(1+\varepsilon)z)$.

Fix a single value of $n$, and let $h,t$ vary.  We claim that the resulting elements $x(n,h,t)$ are all distinct.  Indeed, if $x(n,h,t)=x(n,h',t')$, then rearranging gives
$$\frac{h-h'}{8^{N-n}}=\frac{t-t'}{T}.$$
The coprimality of $8^{N-n},T$ (since $T$ is odd) implies that $T$ divides $t-t'$.  But then in fact $t=t'$, and thence also $h=h'$, as desired.  Next, the number $x(n,h,t)/8^n$ is rational with denominator dividing $8^N \cdot T$.  So different choices of $h,t$ give elements $x(n,h,t)$ that are additively $8^n/(8^N \cdot T)$-separated and hence multiplicatively $(1+\varepsilon)$-separated for the choice $$\varepsilon:=\frac{1}{2 \cdot 8^N \cdot T} \asymp \frac{1}{8^{2N} \cdot N}.$$

We conclude that $M_\varepsilon(A)$ is at least the total number of elements $x(n,h,t)$.  This gives
$$\frac{M_\varepsilon(A)}{|A|} \gg \frac{N \cdot 8^N \cdot (N \cdot 8^N)}{N \cdot 8^{N}}=(N \cdot 8^{2N})^{1/2} \cdot N^{1/2} \asymp \varepsilon^{-1/2}(\log \varepsilon^{-1})^{1/2},$$
which matches the bound in Proposition~\ref{prop:add} up to a constant factor.  Note that $|A|=\varepsilon^{-(1/2+o(1))}$ in this construction.

\section*{Acknowledgments and AI statement}
The author was supported in part by a NSF Mathematical Sciences Postdoctoral Research Fellowship under grant DMS-2501336.

With minimal substantive human input, ChatGPT 6 successively produced upper bounds for $u_n$ of the shape $e^{O(n^{2/3})}$, then $e^{O(\sqrt{n\log n})}$, and finally $e^{O(\sqrt{n})}$.  The argument initially made do with a suboptimal form of the crucial Proposition~\ref{prop:add}, and further prompting produced the sharp form of this proposition.  For this human-written note, I refashioned all of the proofs by hand and added commentary and motivation.

\end{document}